\documentclass[leqno,draft,11pt]{article}
\usepackage{amssymb,amsmath,theorem,a4wide}

\def\dedic{\thanks{Supported by grant IAA100190902 of GA AV \v CR,
project 1M0545 of M\v SMT \v CR,
RVO: 67985840,
and a grant from the John Templeton Foundation.}}

\def\enumup{}

\allowdisplaybreaks[2]

\ifx\loadcyr\undefined\else
\input cyracc.def
\makeatletter
 at 1\@ptsize pt
 at 1\@ptsize pt
 at 1\@ptsize pt
\makeatother

\fi

\def\gobble#1{}
\def\fixsup#1#2{{#1\let\dp\gobble\mathstrut}^#2_}

\def\bme{\hskip.75em\relax}

\let\eq\leftrightarrow
\let\EQ\Leftrightarrow

\def\iff{\quad\text{iff}\quad}
\let\LOR\bigvee
\let\ET\bigwedge
\let\TO\Rightarrow
\def\?{\mathbin?}
\let\ot\leftarrow

\let\model\vDash
\let\nmodel\nvDash

\newbox\circlebox
\setbox\circlebox\hbox{$\bigcirc$}
\def\circled#1{%
  \setbox0\hbox to\wd\circlebox{\hss$#1$\hss}\wd0=0pt
  \box0\copy\circlebox}

\let\fii\varphi
\let\tet\vartheta
\let\ep\varepsilon

\def\greek#1{$\expandafter\greeknum\csname c@#1\endcsname$}
\makeatletter
\def\greeknum#1{\ifcase#1\or\alpha\or\beta\or\gamma\or\delta\or\ep
      \or\digamma\or\zeta\or\eta\or\tet\or\iota\else\@ctrerr\fi}
\makeatother

\def\p#1{\langle#1\rangle}

\def\lh#1{\lvert#1\rvert}

\let\cat\smallfrown

\let\bez\smallsetminus

\let\sset\subseteq
\let\nsset\nsubseteq

\let\onto\twoheadrightarrow

\def\pw#1{\mathcal P(#1)}
\let\nul\varnothing
\def\res{\mathbin\restriction}

\def\twoprimes{\raise.2\fontdimen6\scriptfont2\hbox{$\scriptstyle\prime\prime$}}
\newcommand\rpair[3][3em]{\mathrel{%
   \begin{matrix}%
     \strut\smash{\xrightonto{\hbox to#1{\hss$#2$\hss}}}\\[-1.7ex]%
     \strut\smash{\xleftembed[\hbox to#1{\hss$#3$\hss}]{}}%
   \end{matrix}}}
\makeatletter
\newcommand\xrightonto[2][]{\ext@arrow 0359\rightontofill{#1}{#2}}
\newcommand\xleftembed[2][]{\ext@arrow 3095\leftembedfill{#1}{#2}}
\def\leftembedfill{\arrowfill@\leftarrow\relbar\hookleftnoarrow}
\def\rightontofill{\arrowfill@\relbar\relbar\onto}
\def\hookleftnoarrow{\DOTSB\relbar\joinrel\rhook}
\makeatother

\mathchardef\#="2023 

\ifx\busspf\undefined\else
\usepackage{bussproofs}
\EnableBpAbbreviations

\fi

\ifx\symlasy\undefined

  \def\centerdot#1{%
    \setbox0\hbox{$\mathop{#1}$}\dimen0 \ht0
    \setbox0\hbox{$#1$}\advance\dimen0 -\ht0
    \setbox2\hbox to\wd0{\hss$\mathop{\cdot}$\hss}\wd2=0pt
    \lower\dimen0\box2\box0 }
\else
  
  \def\centerdot#1{%
     \setbox0\hbox{$#1$}%
     \raise0.206\ht0\hbox to\wd0{\hss$\cdot$\hss}%
     \kern-\wd0 \box0 }
\fi

\let\sls|

\def\Up{{\setbox0\hbox{$\uparrow$}%
         \lower\dp0\hbox to\wd0{\hss\vrule width4pt height.4pt\hss}%
         \kern-\wd0\box0}}
\def\UP{{\setbox0\hbox{$\uparrow$}%
         \lower\dp0\hbox to\wd0{\hss\vrule width4pt height.4pt\hss}%
         \kern-\wd0\copy0\kern-\wd0\raise.35ex\box0}}
\def\Down{{\setbox0\hbox{$\downarrow$}%
         \raise\ht0\hbox to\wd0{\hss\vrule width4pt depth.4pt\hss}%
         \kern-\wd0\box0}}

\let\cxt\mathrm

\def\psp{\cxt{PSPACE}}

\def\st{\expandafter\hat}

\let\lgc\mathbf

\def\CPC{\lgc{CPC}}

\def\FL#1{\lgc{FL_{#1}}}

\mathcode`\*="0203
\let\ob\overline

\mathchardef\mhyphen="2D

\def\noproof{\leavevmode\unskip\bme\vadjust{}\nobreak\hfill$\qed$\par}
\let\qed\Box
\newenvironment{Pf}[1][]
  {\par\noindent\textit{Proof\optpar{#1}:}\bme\ignorespaces}
  {\noproof\pagebreak[2]\vskip\medskipamount\ignorespacesafterend}
\def\qedhere{\relax\ifmmode\eqno\qed\expandafter\aftergroup
                   \else\noproof\fi\noqed}
\def\noqed{\let\noproof\relax}

\theoremstyle{plain}
\ifx\shortthm\undefined
\newtheorem{Thm}{Theorem}[section]
\else
\newtheorem{Thm}{Theorem}
\fi

\newtheorem{Lem}[Thm]{Lemma}

\newtheorem{Obs}[Thm]{Observation}

\newtheorem{Cl}{Claim}[Thm]
\def\theCl{\arabic{Cl}}

\theorembodyfont\upshape
\newtheorem{Def}[Thm]{Definition}
\newtheorem{Rem}[Thm]{Remark}
\newtheorem{Exm}[Thm]{Example}

\newenvironment{Pf*}{\let\qed\qedCl\Pf}\endPf

\usepackage[reftex]{theoremref}
\newif\iflinenumbers
\linenumberstrue

{\catcode`\^^I=13 \catcode`\^^M=13
\gdef\doalgo#1#2\end#{\hbox to\hsize{\hss \let^^I\qquad%
  \def\\^^M{\nobreak\hfil\break\vadjust{}\qquad}%
  \fboxsep1em \linenum0 %
  \fbox{\hsize#1\vbox{%
  \everypar{\advance\linenum1 %
      \hbox to1.2em{%
           \hss\iflinenumbers$\scriptstyle\the\linenum$\hskip.6em\fi}}%
  #2}}\hss}\end}}
\newcount\linenum

\def\key{\relax\ifmmode\expandafter\mathbf\else\expandafter\textbf\fi}

\def\allowhyphens{\nobreak\hskip0pt\relax}

\DeclareRobustCommand*\magiclparen{\ifmmode(\else\textup(\allowhyphens\fi}
\DeclareRobustCommand*\magicrparen{\ifmmode)\else\textup)\fi}
\let\lparen=(  \let\rparen=)
\def\magicparon{\catcode`\(\active\catcode`\)\active}
\def\magicparoff{\catcode`\(12 \catcode`\)12 }
\AtBeginDocument{\ifx\ifPreview\iftrue\else\magicparon\fi}
\magicparon
\let (=\magiclparen  \let )=\magicrparen

\ifx\enumup\undefined
  
\else
  
\fi

\def\optpar#1{\ifx\relax#1\relax\else\/ (#1)\fi}

\magicparoff

\mathchardef\comma=\mathcode`\,
{\catcode`\,=\active \gdef,{\comma\penalty\relpenalty}}

\providecommand\dedic{\thanks{Supported by
grant IAA100190902 of GA AV \v CR, Center of Excellence CE-ITI under the grant
P202/12/G061 of GA \v CR, and RVO: 67985840.}}
\author{Emil Je\v r\'abek\dedic\\[\medskipamount]
Institute of Mathematics of the Academy of Sciences\\
\small \v Zitn\'a 25,
115\:67 Praha 1,
Czech Republic,
email: \texttt{jerabek@math.cas.cz}
}

\title{The ubiquity of conservative translations}

\begin{document}
\maketitle
\begin{abstract}
We study the notion of conservative translation between logics
introduced by Feitosa and D'Ottaviano \cite{feit-dott}. We show that
classical propositional logic ($\CPC$) is universal in the sense that
every finitary consequence relation over a countable set of formulas can be
conservatively translated into $\CPC$. The translation is computable
if the consequence relation is decidable. More generally, we show that one
can take instead of $\CPC$ a broad class of logics (extensions of a
certain fragment of full Lambek calculus $\FL{}$) including most
nonclassical logics studied in the literature, hence in a sense,
(almost) any two reasonable deductive systems can be conservatively
translated into each other. We also provide some counterexamples, in
particular the paraconsistent logic $\lgc{LP}$ is not universal.
\end{abstract}

\section{Introduction}\label{sec:introduction}
There have been several proposals of a general concept of a
translation or interpretation between abstract logical systems, see
e.g.\ \cite{cco,mdt} for overviews. A minimalist approach was taken by
da Silva, D'Ottaviano and Sette \cite{sos}: a logic (deductive system)
is given by any Tarski-style consequence operator, and then a
translation of one logic in another is an arbitrary mapping of
formulas to formulas preserving the consequence relation. Feitosa and
D'Ottaviano \cite{feit-dott} consider the stronger notion of
\emph{conservative translations}, which preserve the consequence
relation in both directions. This avoids uninteresting examples of
translations such as mapping all formulas to a fixed tautology.

This notion of a conservative translation is still very general
(perhaps too general): for instance, translations are not required
to respect the structure of formulas in any way, to be computable, or
to preserve any properties of the logic.
For this reason, it is natural to expect that there should exist a
conservative translation between more or less any two
reasonable deductive systems. Nevertheless, no result to such effect
appears in the literature. Instead, there are several
papers presenting proofs (often non-constructive) of the existence of
conservative translations between particular pairs of logics:
\cite{dott-feit:mv,dott-feit:para,dott-feit:luk,dott-feit:persp}.

The main purpose of this paper is to demonstrate that for a quite
large class of logics, it is indeed possible to construct a
conservative translation between any two of them. First, we prove that
an arbitrary finitary deductive system $L$ in countably many formulas can be
conservatively translated into classical propositional logic ($\CPC$),
in either the single-conclusion or multiple-conclusion setting.
Our translation is constructed by an explicit inductive definition,
and it is computable whenever $L$ is decidable. Moreover, the
translation has the additional property of being most general in the
sense that every other translation of $L$ to $\CPC$ is equivalent to its
substitution instance.

Let us define a logic $L$ to be \emph{universal} if every finitary
deductive system in countably many formulas can be conservatively
translated into $L$, so that the result above can be restated by
saying that $\CPC$ is universal. We generalize this result by showing
that every deductive system between the $\to,\ot,\land$ fragment of the full
Lambek calculus $\FL{}$ (see \cite{reslat}) and the corresponding
fragment of $\CPC$ is universal, and similarly, any deductive system between
$\lgc{BCK}$ (the implication fragment of $\FL{ew}$) and $\CPC\res_\to$
is universal. This establishes the universality of most of
nonclassical logics studied in the literature, as they typically
extend (a suitable fragment of) $\FL{}$ in one way or another: this
includes e.g.\ intuitionistic and intermediate logics, various modal,
substructural, fuzzy, or relevant logics, both propositional and
first-order.

As an additional example, we show that Kleene's logic with truth
constants is universal, whereas the paraconsistent logic
$\lgc{LP}$---based on the same algebra but with a different choice of
designated truth values---is \emph{not} universal. We also completely
characterize universal fragments of $\CPC$: a fragment $\CPC\res_B$ is
universal if and only if implication is definable from $B$. (In
particular, we obtain a couple of nontrivial examples of logics into
which $\CPC$ cannot be conservatively translated, namely $\lgc{LP}$
and the fragments $\CPC\res_{\eq,\neg}$,
$\CPC\res_{\land,\lor,\bot,\top}$.)

The paper is organized as follows. In Section~\ref{sec:class-logic} we
give basic definitions, and we construct conservative
translations into classical logic. In
Section~\ref{sec:universal-logics} we investigate the class of
universal deductive systems, as detailed above. Section~\ref{sec:conclusion}
consists of concluding remarks.

\section{Translation to classical logic}\label{sec:class-logic}
In this section, we are going to construct conservative translations
of (almost) arbitrary logics into classical logic
(\th\ref{thm:cpc-1}). First we review the relevant definitions to fix
the notation.
\begin{Def}
A pair $L=\p{F,\vdash}$ is a \emph{deductive system} (or \emph{logic})
over a set of \emph{formulas} $F$ if ${\vdash}\sset\pw F\times F$ is a
Tarski-style \emph{consequence relation}, i.e., if it satisfies
\begin{enumerate}
\item $\fii\vdash\fii$,
\item $\Gamma\vdash\fii$ implies $\Gamma,\Gamma'\vdash\fii$,
\item if $\Gamma\vdash\fii$ and $\Delta\vdash\psi$ for all
$\psi\in\Gamma$, then $\Delta\vdash\fii$,
\end{enumerate}
for every $\fii\in F$ and $\Gamma,\Gamma',\Delta\sset F$. A deductive system is
\emph{finitary} if
\begin{enumerate}
\setcounter{enumi}3
\item $\Gamma\vdash\fii$ implies $\Gamma'\vdash\fii$ for some finite
$\Gamma'\sset\Gamma$.
\end{enumerate}
When discussing algorithmic issues, we will tacitly assume that $F$ is
encoded as a recursively enumerable subset of $\omega$. If
$X=\p{\fii_1,\dots,\fii_n}$ is a sequence of formulas, we will also
write $X\vdash\fii$ instead of $\{\fii_1,\dots,\fii_n\}\vdash\fii$ by
abuse of notation (we will never use $\vdash$ for a sequent arrow).

A deductive system $L=\p{F,\vdash}$ is a \emph{propositional logic} if
$F$ is the set of formulas built inductively from a set of variables
and a set of finitary connectives (i.e., $F$ is a free algebra in a
particular signature), and $\vdash$ is structural (substitution-invariant):
\begin{enumerate}
\setcounter{enumi}4
\item $\Gamma\vdash\fii$ implies $\sigma(\Gamma)\vdash\sigma(\fii)$
for every substitution $\sigma$,
\end{enumerate}
where as usual, a substitution is a homomorphism of free algebras.
Let $\CPC=\p{F_\CPC,\vdash_\CPC}$ denote the usual consequence
relation of classical propositional logic in countably infinitely many
variables using an arbitrary functionally complete finite set of Boolean
connectives (the exact choice does not matter, as classical
consequence is unaffected by translation of formulas to a language
with a different set of basic connectives). In $\CPC$, we will employ big
conjunctions and disjunctions $\ET\Gamma$, $\LOR\Gamma$ in the usual
way; in particular, $\ET\nul=\top$ and $\LOR\nul=\bot$ (again, it does
not matter whether these constants are included in the set of basic
connectives, or defined by equivalent more complicated formulas, even if they
involve extra variables).
\end{Def}
\begin{Def}
A \emph{translation} from a deductive system $L_0=\p{F_0,\vdash_0}$ to
a deductive system
$L_1=\p{F_1,\vdash_1}$ is a function $f\colon F_0\to F_1$ such that
$$\Gamma\vdash_0\fii\ \TO\ f(\Gamma)\vdash_1f(\fii)$$
for every $\Gamma\sset F_0$, $\fii\in F_0$. We will write this as
$f\colon L_0\to L_1$. The translation $f$ is
\emph{conservative}, written as $f\colon L_0\to_cL_1$, if
$$\Gamma\vdash_0\fii\ \EQ\ f(\Gamma)\vdash_1f(\fii).$$
We write $L_0\le_cL_1$ if there exists a conservative translation
$f\colon L_0\to_cL_1$.
\end{Def}

The translations to classical logic we construct have an additional
property which might be of independent interest, hence we give it a name:
\begin{Def}
Let $L_0=\p{F_0,\vdash_0}$ be a deductive system, and
$L_1=\p{F_1,\vdash_1}$ a propositional logic. A translation $f\colon
L_0\to L_1$ is \emph{most general} if for every translation $g\colon
L_0\to L_1$, there exists a substitution $\sigma$ such that
$g(\fii)\dashv\vdash_1\sigma(f(\fii))$ for every $\fii\in F_0$.
\end{Def}
Notice that if $L_0\le_cL_1$ and $L_1$
is finitary, then $L_0$ is also finitary.

The main result of this section is:
\begin{Thm}\th\label{thm:cpc-1}
For every finitary deductive system $L=\p{F,{\vdash}}$ over a countable set of
formulas $F$, there exists a conservative most general translation
$f\colon L\to_c\CPC$.

If\/ $\vdash$ is decidable, then $f$ is computable. In general, $f$
is Turing equivalent to (the finitary fragment of) $\vdash$.
\end{Thm}

We will prove \th\ref{thm:cpc-1} below as a corollary to its
multiple-conclusion version. Apart from being more
general, the construction of the translations in the
multiple-conclusion case is more transparent and displays better the
underlying symmetry, we thus find it preferable to giving a direct
proof for the single-conclusion case, which feels a bit ad hoc.
\begin{Def}\th\label{def:mc}
A pair $L=\p{F,\vdash}$ is a \emph{multiple-conclusion deductive
system} (or \emph{multiple-conclusion logic})
\cite{sh-sm} if ${\vdash}\sset\pw F\times\pw F$ satisfies
\begin{enumerate}
\item $\fii\vdash\fii$,
\item $\Gamma\vdash\Delta$ implies $\Gamma,\Gamma'\vdash\Delta,\Delta'$,
\item\label{item:cut} if $\Gamma,\Pi\vdash\Lambda,\Delta$ for every
$\Pi,\Lambda$ such
that $\Pi\cup\Lambda=\Xi$, then $\Gamma\vdash\Delta$,
\end{enumerate}
for every $\fii\in F$ and $\Gamma,\Gamma',\Delta,\Delta',\Xi\sset F$.
(Condition \eqref{item:cut} is a form of the cut rule.)
A multiple-conclusion deductive system is
\emph{finitary} if
\begin{enumerate}
\setcounter{enumi}3
\item $\Gamma\vdash\Delta$ implies $\Gamma'\vdash\Delta'$ for some finite
$\Gamma'\sset\Gamma$, $\Delta'\sset\Delta$.
\end{enumerate}
Note that if $L$ is finitary, condition \eqref{item:cut} can be
equivalently simplified to
\begin{enumerate}
\item[(\ref{item:cut}$'$)] if $\Gamma,\fii\vdash\Delta$ and
$\Gamma\vdash\fii,\Delta$, then $\Gamma\vdash\Delta$.
\end{enumerate}
$L$ is \emph{consistent} if $\nul\nvdash\nul$.
Let $\CPC_m=\p{F_\CPC,\model_\CPC}$ denote the maximal structural
multiple-conclusion consequence 
relation for classical propositional logic: $\Gamma\model_\CPC\Delta$ iff there is no
$0$--$1$ assignment $v$ such that $v(\fii)=1$ for all $\fii\in\Gamma$ and
$v(\psi)=0$ for all $\psi\in\Delta$. (In other words,
$\Gamma\model_\CPC\Delta$ iff there are finite subsets
$\Gamma'\sset\Gamma$, $\Delta'\sset\Delta$ such that
$\vdash_\CPC\ET\Gamma'\to\LOR\Delta'$.)

We generalize the notions of translations, conservative translations,
propositional logics, and most general translations to the
multiple-conclusion setting in the obvious way.
\end{Def}
\begin{Thm}\th\label{thm:cpc-m}
For every finitary consistent multiple-conclusion deductive system
$L=\p{F,{\vdash}}$ over a countable set of formulas $F$, there exists a
conservative most general translation $f\colon L\to_c\CPC_m$.

If\/ $\vdash$ is decidable, then $f$ is computable. In general, $f$
is Turing equivalent to $\vdash$.
\end{Thm}
\begin{Pf}
Let $F=\{\alpha_n:n\in\omega\}$ be a (not necessarily injective)
enumeration. We will define a sequence of formulas $f(\alpha_n):=\beta_n\in
F_\CPC$ by induction on $n$. We denote by $p_n$ the $n$th
propositional variable of $\CPC$. We abbreviate
$\ob\alpha_X=\{\alpha_i:i\in X\}$, and similarly for $\ob\beta_X$; moreover,
we are going to use the identity $n=\{i\in\omega:i<n\}$.

Assume by the induction hypothesis that $\beta_i$ have been already
defined for all $i<n$ in such a way that
\begin{equation}\label{eq:1}
\ob\alpha_X\vdash\ob\alpha_Y\ \TO\ \ob\beta_X\model_\CPC\ob\beta_Y
\end{equation}
for every $X,Y\sset n$. (This holds for $n=0$ as $\vdash$ is
consistent by assumption.) Define\footnote{Note that the $\beta_n$ are
defined by complete (ordinal) induction, which requires no separate
base case. For example, the given definition implies that $\beta_0$ is
one of $\top\lor p_0\land\top$ ($\equiv\top$), $\bot\lor p_0\land\bot$
($\equiv\bot$), or $\bot\lor p_0\land\top$ ($\equiv p_0$), depending on whether
$\nul\vdash\alpha_0$, $\alpha_0\vdash\nul$, or neither, respectively.}
$\beta_n$ by
\begin{align*}
\gamma_n&:=\LOR_{
  \substack{X,Y\sset n\\\ob\alpha_X\vdash\alpha_n,\ob\alpha_Y}}
 \left(\ET\ob\beta_X\land\neg\LOR\ob\beta_Y\right),\\
\delta_n&:=\ET_{
  \substack{X,Y\sset n\\\ob\alpha_X,\alpha_n\vdash\ob\alpha_Y}}
 \left(\ET\ob\beta_X\to\LOR\ob\beta_Y\right),\\
\beta_n&:=\gamma_n\lor p_n\land\delta_n.
\end{align*}
Notice that it does not matter whether we read the definition
of $\beta_n$ as $\gamma_n\lor(p_n\land\delta_n)$ or $(\gamma_n\lor
p_n)\land\delta_n$, as $\model_\CPC\gamma_n\to\delta_n$: if
$\ob\alpha_X\vdash\alpha_n,\ob\alpha_Y$ and
$\ob\alpha_W,\alpha_n\vdash\ob\alpha_Z$, then
$\ob\alpha_X,\ob\alpha_W\vdash\ob\alpha_Y,\ob\alpha_Z$, hence
$$\model_\CPC\ET\ob\beta_X\land\ET\ob\beta_W\to\LOR\ob\beta_Y\lor\LOR\ob\beta_Z$$
by \eqref{eq:1}, i.e.,
$$\model_\CPC\left(\ET\ob\beta_X\land\neg\LOR\ob\beta_Y\right)
   \to\left(\ET\ob\beta_W\to\LOR\ob\beta_Z\right).$$

We claim that \eqref{eq:1} holds for $X,Y\sset n+1$. If $X,Y\sset n$,
this follows from the induction hypothesis. If $n\in X\cap Y$, then
trivially $\ob\beta_X\model_\CPC\ob\beta_Y$. Assume that
$\ob\alpha_X\vdash\alpha_n,\ob\alpha_Y$, where $X,Y\sset n$. Then
$$\ET\ob\beta_X\land\neg\LOR\ob\beta_Y\model_\CPC\gamma_n\model_\CPC\beta_n,$$
hence by reasoning in $\CPC$,
$$\ob\beta_X\model_\CPC\ob\beta_Y,\beta_n.$$
The case $\ob\alpha_X,\alpha_n\vdash\ob\alpha_Y$ is handled similarly using
the definition of $\delta_n$.

Thus, $f$ is well defined, and by \eqref{eq:1} and finitarity of $L$,
it is a translation of $L$ to $\CPC_m$. In order to show that $f$ is
conservative, assume that $\ob\alpha_W\nvdash\ob\alpha_Z$, we need to prove
$\ob\beta_W\nmodel_\CPC\ob\beta_Z$. Obviously, $W\cap Z=\nul$. By the cut
rule (i.e., \th\ref{def:mc} \eqref{item:cut}, applied with $\Xi=F$),
we may assume that $W\cup Z=\omega$. Let $v$ be the valuation such
that
$$v(p_n)=\begin{cases}1,&n\in W,\\0,&n\in Z.\end{cases}$$
We will show $v(\beta_n)=v(p_n)$ by induction on $n$, which implies
$\ob\beta_W\nmodel_\CPC\ob\beta_Z$.

Assume that $n\in W$. If $X,Y\sset n$ are such that
$\ob\alpha_X,\alpha_n\vdash\ob\alpha_Y$, then we cannot have simultaneously
$X\sset W$ and $Y\sset Z$. If $i\in X\bez W$, then $v(\beta_i)=0$ by
the induction hypothesis; similarly, if $i\in Y\bez Z$, then
$v(\beta_i)=1$. Thus, $v(\ET\ob\beta_X)=0$ or $v(\LOR\ob\beta_Y)=1$.
Since $X,Y$ were arbitrary, we obtain $v(\beta_n)=v(\delta_n)=1$.

If $n\in Z$, we obtain $v(\beta_n)=v(\gamma_n)=0$ by a similar argument.

Clearly, the explicit recursive definition of $f$ can be realized by an
algorithm with an oracle for $\vdash$. On the other hand, since $f$ is
a conservative translation into the decidable logic $\CPC_m$, the
relation $\vdash$
is Turing reducible to $f$ (or its graph, if we insist on oracles
being sets rather than functions).

It remains to show that $f$ is a most general translation of $L$ to
$\CPC_m$. Let $g\colon L\to\CPC_m$, and let $\sigma$ be the
substitution defined by $\sigma(p_n)=g(\alpha_n)$. We have to establish
\begin{equation}\label{eq:2}
\model_\CPC g(\alpha_n)\eq\sigma(\beta_n)
\end{equation}
for every $n$, and we proceed by induction on $n$. The definition of
$\beta_n$ implies that \eqref{eq:2} is equivalent to
$$\sigma(\gamma_n)\model_\CPC g(\alpha_n)\model_\CPC\sigma(\delta_n).$$
Using the definitions of $\gamma_n,\delta_n$ and the induction
hypothesis, this is equivalent to
$$\LOR_{
  \substack{X,Y\sset n\\\ob\alpha_X\vdash\alpha_n,\ob\alpha_Y}}
 \left(\ET g(\ob\alpha_X)\land\neg\LOR g(\ob\alpha_Y)\right)
\model_\CPC g(\alpha_n)\model_\CPC
\ET_{
  \substack{X,Y\sset n\\\ob\alpha_X,\alpha_n\vdash\ob\alpha_Y}}
 \left(\ET g(\ob\alpha_X)\to\LOR g(\ob\alpha_Y)\right),$$
which in turn follows from the fact that $g$ is a translation: for
example, if $X,Y\sset n$ are such that
$\ob\alpha_X\vdash\alpha_n,\ob\alpha_Y$, then $g(\ob\alpha_X)\model_\CPC
g(\alpha_n),g(\ob\alpha_Y)$, hence $\ET g(\ob\alpha_X)\land\neg\LOR
g(\ob\alpha_Y)\model_\CPC g(\alpha_n)$.
\end{Pf}
\begin{Pf}[of \th\ref{thm:cpc-1}]
Let $L=\p{F,\vdash}$ be a finitary deductive system over countable $F$, and
define its conservative multiple-conclusion extension $L^m=\p{F,{\vdash}^m}$ by
$$\Gamma\vdash^m\Delta\iff\exists\psi\in\Delta\:\Gamma\vdash\psi.$$
Let $f\colon L^m\to_c\CPC^m$ be the conservative minimal translation
from \th\ref{thm:cpc-m}. Since $\vdash_\CPC$ is the single-conclusion
fragment of $\model_\CPC$, $f\colon L\to_c\CPC$.
Moreover, if $g\colon L\to\CPC$, then $g\colon L^m\to\CPC_m$, hence
$g$ is $\CPC$-equivalent to $\sigma\circ f$ for some substitution $\sigma$.
\end{Pf}
\begin{Rem}
Even if $L=\CPC$, the translation $f$ from \th\ref{thm:cpc-1} is not
(equivalent to) the identity, since it has the additional property
that $f(\Gamma)\vdash_\CPC\LOR_{i<n}f(\fii_i)$ implies
$f(\Gamma)\vdash_\CPC f(\fii_i)$ for some $i<n$.
\end{Rem}
\begin{Rem}
Let us estimate the complexity of our translation $f$. Assume that
formulas are represented by strings in a finite alphabet and
enumerated in the natural way so that shorter formulas
have smaller index, and consider a formula $\fii$ of length $n$, so
that $\fii=\alpha_m$ for some $m=2^{O(n)}$. In order to compute $f(\fii)$,
we have to determine whether $\Gamma\vdash\Delta$ for sets
$\Gamma,\Delta$ of formulas whose length is at most $n$; there are
$2^{O(n)}$ such formulas, hence $2^{2^{O(n)}}$ sets. Thus, we define
$f(\fii)$ by iteration of length $2^{O(n)}$, and in each step, the
formula $\beta_i$ is constructed from at most $2^{2^{O(n)}}$ copies of
formulas constructed earlier. It follows that the total length of
$f(\fii)$ is $(2^{2^{O(n)}})^{2^{O(n)}}=2^{2^{O(n)}}$.

Assume that $\vdash$ is decidable in
$\cxt{EXP}=\cxt{DTIME}(2^{n^{O(1)}})$. Since each of the sequents
$\Gamma\vdash\Delta$ above has size $s=2^{O(n)}$, $f(\fii)$ is
computable in time $2^{2^{O(n)}}2^{s^{O(1)}}=2^{2^{O(n)}}$. Moreover,
if we compute $f$ in a left-to-right fashion by a recursive procedure
mimicking its definition, we need recursion depth $2^{O(n)}$ and local
storage $2^{O(n)}$ for each recursive call, plus the space needed to
check $\Gamma\vdash\Delta$. Thus, if ${\vdash}\in\psp$, then
$f(\fii)$ is computable in space $2^{O(n)}$.

In general, if $\vdash$ is computable in time $t(n)$ and space
$s(n)$, where both $t$ and $s$ are monotone, then $f(\fii)$ is
computable in time $2^{2^{O(n)}}t(2^{O(n)})$ and space $2^{O(n)}+s(2^{O(n)})$.
\end{Rem}

\section{Universal logics}\label{sec:universal-logics}
From now on, we only consider single-conclusion logics.
\begin{Def}
A deductive system $L_0$ is \emph{universal}
if $L\le_cL_0$ for every finitary deductive system $L$ over countably
many formulas.

(In all cases where we establish universality below, it is possible to
construct an $f\colon L\to_cL_0$ Turing
equivalent to $\vdash_L$, as in
the case of $\CPC$. However, we decided not to include this condition in
the definition.)
\end{Def}
In the previous section, we proved that $\CPC$ is a universal logic.
The main result of this section, \th\ref{thm:univ}, is a
generalization of our construction to a large class of nonclassical
logics in place of $\CPC$.

Before we get to the main result, we discuss some examples showing
that the question of which logics $\CPC$ can be conservatively translated
into is considerably more subtle than which logics can be conservatively
translated into $\CPC$.

From \th\ref{thm:cpc-1} and the transitivity of $\le_c$, we immediately obtain:
\begin{Obs}\
\begin{enumerate}
\item If $L_0$ is universal and $L_0\le_cL_1$, then $L_1$ is universal.
\item $L_0$ is universal iff\/ $\CPC\le_cL_0$.
\end{enumerate}
\end{Obs}
\begin{Exm}\th\label{exm:kleene}
Let $A_3=\p{\{0,*,1\},\land,\lor,0,1,\neg}$ be the $3$-element bounded
lattice endowed with an operation $\neg0=1$, $\neg1=0$, $\neg*=*$.
Recall that \emph{Kleene's $3$-valued logic $\lgc K$} (with truth constants) is
the propositional logic using connectives $\land,\lor,\neg,\top,\bot$
whose consequence relation is defined by $A_3$ where $1$ is the
only designated value, and the paraconsistent \emph{logic of paradox
$\lgc{LP}$} is defined similarly but with both $1,*$ taken as
designated. We have:
\begin{enumerate}
\item\label{item:k}
$\lgc K$ is universal.
\item\label{item:lp}
$\lgc{LP}$ is not universal.
\end{enumerate}
(Notice that $\lgc K$ without truth constants is trivially not
universal, as it has no tautologies.)
\end{Exm}
\begin{Pf}
\eqref{item:k}: Let $f(\fii)$ be a conjunctive normal form of $\fii$,
obeying the convention that no variable and its negation can appear
simultaneously in a clause, and we use $\top,\bot$ for empty
conjunctions and disjunctions, respectively. We claim that
$f\colon\CPC\to_c\lgc K$. Since $\lgc K\sset\CPC$ and
$\fii\dashv\vdash_\CPC f(\fii)$, we clearly have
$$f(\Gamma)\vdash_{\lgc K}f(\fii)\ \TO\ \Gamma\vdash_\CPC\fii.$$
In order to show the converse implication, it suffices to prove that
$$\Gamma\vdash_\CPC\fii\ \TO\ \Gamma\vdash_{\lgc K}\fii$$
holds whenever $\Gamma\cup\{\fii\}$ is a set of clauses.
Let $v$ be a valuation in $A_3$ such that
$v(\Gamma)=1$ and $v(\fii)\ne1$. We modify $v$ to make it a Boolean
valuation $v'$ as follows. If $l\in\fii$ is a literal such that $v(l)=*$,
we put $v'(l)=0$; we can do this for all such $l$
simultaneously since $\fii$ does not contain both $l$ and $\neg l$. If
$p$ is a variable such that $v(p)=*$ and neither $p$ nor $\neg p$
appears in $\fii$, we pick $v'(p)\in\{0,1\}$ arbitrarily. After this
modification, $v'$ is a classical valuation such that $v'(\fii)=0$, and
since all literals with value $1$ kept their value, we still have
$v'(\Gamma)=1$.

\eqref{item:lp}:
Assume for contradiction $f\colon\CPC\to_c\lgc{LP}$. Let
$\{v_i:i<n\}$ be the list of all
valuations in $A_3$ such that $v_i(f(\bot))=0$ with $v_i(p_j)=*$ for
every variable $p_j$ not occurring in $f(\bot)$.
Put $\fii_i=p_i$ for $i<n$,
$\fii_n=\neg\ET_{i<n}p_i$. We have
$\fii_0,\dots,\fii_n\vdash_\CPC\bot$,
which implies $f(\fii_0),\dots,f(\fii_n)\vdash_\lgc{LP}f(\bot)$. Since
$v_i(f(\bot))=0$, we must have $v_i(f(\fii_{j_i}))=0$ for some $j_i\le
n$. Put $J=\{j_i:i<n\}$. We claim that
$$\{f(\fii_j):j\in J\}\vdash_\lgc{LP}f(\bot).$$
Indeed, if $v(f(\bot))=0$, there exists an $i$ such that $v$ and
$v_i$ coincide on variables occurring in $f(\bot)$. We have
$v_i(f(\fii_j))=0$ for some $j\in J$. If $\prec$ is the partial order
induced by $*\prec0$, $*\prec1$, then functions definable in $A_3$ are
$\prec$-monotone, and $v_i\preceq v$, hence $v(f(\fii_j))=0$.
Thus, by the conservativity of $f$,
$$\{\fii_j:j\in J\}\vdash_\CPC\bot\text{ for some }\lh J\le n.$$
This contradicts the definition of $\fii_0,\dots,\fii_n$.
\end{Pf}

In order to get some insight which logics can or cannot be expected to
be universal depending on their available list of connectives, we
characterize universal fragments of classical logic below.
\begin{Def}
If $L$ is a propositional logic, and $B$ a set of connectives
definable in $L$ (i.e., $L$-formulas), we denote by $L\res_B$ the
fragment of $L$ using only formulas built from $B$. (We treat
$L\res_B$ as having $B$ as the basic set of connectives, regardless of
the basic connectives of~$L$.)

A \emph{clone} on a set $X$ is a set of finitary
operations on $X$ which is closed under composition and contains all
projections. If $B$ is a set of operations on $X$, then we denote by
$[B]$ the clone generated by $B$. Notice that if $B$ is a set of
Boolean functions, then definable functions in $\CPC\res_B$ are
exactly the functions from $[B]$, hence clones on $\{0,1\}$
are in $1$--$1$ correspondence with fragments of $\CPC$ considered up to term
equivalence.

The lattice of clones on $\{0,1\}$ was completely described by Post
\cite{post} (see also Lau \cite{lau} for a modern exposition). We will
in particular need to refer to the following clones:
\begin{itemize}
\item The clone $P_0$ of all $0$-preserving functions (i.e., $f(0,\dots,0)=0$).
\item The clone $D$ of all self-dual functions (i.e., $f(\neg x_1,\dots,\neg
x_n)=\neg f(x_1,\dots,x_n)$).
\item The clone $A$ of all affine functions ($f(x_1,\dots,x_n)=\sum_{i\in
I}x_i+c$, where $c\in\{0,1\}$, $I\sset\{1,\dots,n\}$, and $+$ denotes
addition modulo $2$).
\item The clone $M$ of all monotone functions.
\item The clone $T_1^\infty$ of functions bounded below by a variable
(there exists $i$ such that $x_i\le f(x_1,\dots,x_n)$ for every $\vec
x\in\{0,1\}^n$).
\end{itemize}
\end{Def}
The following lemma follows immediately from inspection of Post's
lattice, though we invite the reader to give a direct proof:
\begin{Lem}\th\label{lem:post}
If $B$ is a set of Boolean functions, then ${\to}\notin[B]$ if and only
if $B$ is included in $P_0$, $D$, $A$, or $M$.
\noproof\end{Lem}
\begin{Thm}\th\label{thm:cpc-frag}
Let $B$ be a set of Boolean functions. The fragment $\CPC\res_B$ is universal
if and only if $\to$ is definable from $B$.
\end{Thm}
\begin{Pf}
Left-to-right: if ${\to}\notin[B]$, then $B$ is included in one of the
clones mentioned in \th\ref{lem:post}. If $B\sset P_0$ or $B\sset D$,
then $\top\notin[B]$; in other words, $\CPC\res_B$ has no tautologies,
and therefore cannot be universal.

Let $B\sset A$. We claim that if
$\fii,\psi\in A$, then $\fii\le\psi$ only if $\fii=0$
or $\psi=1$ or $\fii=\psi$. Write $\fii(\vec x)=\sum_{i\in I}x_i+c$,
$\psi(\vec x)=\sum_{i\in J}x_i+d$. If $\psi\ne1$, there is a Boolean
valuation $v$ such that $v(\psi)=0$. If $I\nsset J$, we can change the
valuation of any $x_i$ such that $i\in I\bez J$ to make $v(\fii)=1$,
contradicting $\fii\le\psi$. Thus, $\psi=1$ or $I\sset J$. Since
$\fii\le\psi$ implies $\neg\psi\le\neg\fii$, the same argument gives
$\fii=0$ or $J\sset I$. Finally, if $I=J$, then $\psi=\fii$ or
$\psi=\neg\fii$; in the latter case, $\fii\le\psi$ can only hold if
$\fii$ and $\psi$ are constant functions $0$ and $1$, respectively.

In particular, there is no strictly
increasing chain of length more than $3$ of affine functions ordered
by entailment, hence we cannot conservatively translate $\CPC$ (even
with just $2$ variables) into $\CPC\res_B$.

If $B\sset M$, we will show $\CPC\res_B\le_c\lgc{LP}$, hence
$\CPC\res_B$ is not universal by \th\ref{exm:kleene}. Since
$M=[\land,\lor,\top,\bot]$, we may assume
$B=\{\land,\lor,\top,\bot\}$. Let $\sigma$ be the substitution such
that
$\sigma(p)=p\land\neg p$. We claim
$$\Gamma\vdash_{\CPC\res_B}\fii\ \EQ\ 
  \sigma(\Gamma)\vdash_\lgc{LP}\sigma(\fii).$$
Notice that $\sigma$, being a substitution, is a bounded lattice
homomorphism of the
respective free algebras. Let $A_2$ denote the $2$-element bounded
lattice. The mapping $\pi\colon A_3\to A_2$ such that
$\pi(1)=\pi(*)=1$, $\pi(0)=0$, is also a bounded lattice homomorphism,
and it preserves (in both directions) the sets of designated elements.

If $v$ is a valuation in $A_3$ such that $v(\sigma(\Gamma))\ge*$,
$v(\sigma(\fii))=0$, then $v'=\pi\circ v\circ\sigma$ is a valuation in
$A_2$ such that $v'(\Gamma)=1$ and $v'(\fii)=0$, hence
$\Gamma\nvdash_\CPC\fii$.

Conversely, if $v'$ is a valuation in $A_2$ such that $v'(\Gamma)=1$
and $v'(\fii)=0$, let $v$ be the valuation in $A_3$ induced by
$$v(p_i)=\begin{cases}*&v'(p_i)=1,\\0&v'(p_i)=0.\end{cases}$$
Then $\pi\circ v\circ\sigma=v'$, hence $v(\sigma(\Gamma))\ge*$,
$v(\sigma(\fii))=0$.

Right-to-left: we construct $f\colon\CPC\to\CPC\res_\to$ as follows.
First, we rename all propositional variables in the style of Hilbert's
hotel so that we obtain a spare variable $q$ which does not occur in
any formulas. Then, for each formula $\fii$ not containing $q$, let
$f(\fii)$ be an implicational formula equivalent to $\fii\lor q$; it
exists as $[\to]=T_1^\infty$. (For a more explicit construction, we can
use the functional completeness of $\{\to,\bot\}$ to write $\fii(\vec
p)\eq\psi(\vec p,\bot)$ for some $\psi\in[\to]$, and then put
$f(\fii)=(\psi(\vec p,q)\to q)\to q$.) It is easy to see that
$f\colon\CPC\to_c\CPC\res_\to$.
\end{Pf}
\begin{Rem}
Ideally, we would like to prove that a logic is universal whenever it
meets some simple general conditions, such as those studied in
abstract algebraic logic (see \cite{czel}).
However, on the one hand, the affine fragments $\CPC\res_\eq$ or
$\CPC\res_{\eq,\neg}$ are strongly regularly finitely
algebraizable (i.e., as nice as it can get from the point of
view of AAL), on the other hand, Kleene's logic is not even equivalential.
This shows that universality does not have much to do with abstract
algebraic properties of the logic. Consequently, if we want to establish
universality of a class of logics, we cannot rely only on their
general properties, at some point we have to resort to working with
particular systems. We at least try to pick as weak a base system as
possible so that our result covers a broad class of logics including
most systems studied in the literature.
\end{Rem}
\begin{Def}\th\label{def:fl}
A \emph{residuated lattice} is a structure
$\p{L,\land,\lor,\cdot,\to,\ot,1}$ where $\p{L,\land,\lor}$ is a
lattice, $\p{L,\cdot,1}$ is a monoid, and
$$b\le a\to c\ \EQ\ a\cdot b\le c\ \EQ\ a\le c\ot b$$
for every $a,b,c\in L$. (In particular, $a\cdot(a\to b)\le b$, $(b\ot
a)\cdot a\le b$). An \emph{FL-algebra} is a residuated lattice $L$
with a distinguished point $0\in L$. The \emph{full Lambek calculus
$\FL{}$} is the propositional logic using connectives
$\land,\lor,\cdot,\to,\ot,1,0$ such that $\vdash_\FL{}$ is complete
with respect to the class of logical matrices whose underlying
algebras are FL-algebras $L$, with $\{x\in L:x\ge1\}$ taken as the set
of designated elements. $\FL e$ is complete with respect to
commutative FL-algebras ($x\cdot y=y\cdot x$), and $\FL{ew}$ with
respect to $0$-bounded integral ($0\le x\le 1$) commutative FL-algebras.
In a sequent calculus formulation of $\FL{}$, $\lgc e$ corresponds to
the exchange rule, and $\lgc w$ to the weakening rule.
For more information about $\FL{}$ and its extensions or fragments, we
refer the reader to \cite{reslat}.

If $\Gamma=\p{\fii_1,\dots,\fii_k}$ is a sequence of formulas, we define
\begin{align*}
\Gamma\to\psi&:=\fii_1\to(\fii_2\to(\fii_3\to\cdots(\fii_k\to\psi)\cdots)),\\
\psi\ot\Gamma&:=(\cdots((\psi\ot\fii_1)\ot\fii_2)\cdots\ot\fii_{k-1})\ot\fii_k.
\end{align*}
If $k=0$, it is understood that $\Gamma\to\psi=\psi\ot\Gamma=\psi$. We
also put $\prod\Gamma=\fii_1\cdot\fii_2\cdot\ldots\cdot\fii_k$ ($\prod\Gamma=1$
if $k=0$), and
$\fii^k=\underbrace{\fii\cdot\fii\cdot\ldots\cdot\fii}_{k\text{ times}}$.
Notice that in $\FL{}$, $\Gamma\to\psi$ is equivalent (in the sense of
obtaining the same value under any valuation in any $\FL{}$-algebra) to
$\prod\Gamma^{-1}\to\psi$, and $\psi\ot\Gamma$ is equivalent to
$\psi\ot\prod\Gamma^{-1}$, where $\Gamma^{-1}$ denotes the reversal of
the sequence $\Gamma$.
\end{Def}
\begin{Rem}\th\label{rem:transl}
Let $L$ be a finitary deductive system over countably many formulas
$F=\{\alpha_n:n<\omega\}$. From the proofs of
\th\ref{thm:cpc-1,thm:cpc-m} we know that there is
a conservative translation $f\colon L\to_c\CPC$ such that
$f(\alpha_n)=\beta_n$ is inductively defined to be equivalent to the formula
\begin{equation}\label{eq:tr}
\ET_{\substack{X\sset n>k\\\ob\alpha_X,\alpha_n\vdash_L\alpha_k}}
         (\ob\beta_X\to\beta_k)
\land\biggl(p_n\lor
\LOR_{\substack{Z\sset n\\\ob\alpha_Z\vdash_L\alpha_n}}\ET\ob\beta_Z\biggr).
\end{equation}
\end{Rem}
\begin{Thm}\th\label{thm:univ}
A deductive system is universal whenever it conservatively extends a
deductive system
$L_0$ such that
\begin{enumerate}
\item\label{item:fl} $\FL{}\res_{\to,\ot,\land}\sset L_0\sset\CPC\res_{\to,\ot,\land}$, or
\item\label{item:fle} $\FL e\res_{\to,\land}\sset L_0\sset\CPC\res_{\to,\land}$, or
\item\label{item:bck} $\FL{ew}\res_\to=\lgc{BCK}\sset L_0\sset\CPC\res_\to$.
\end{enumerate}
\end{Thm}
\begin{Pf}
\eqref{item:fl}:
Let $L$ be as in \th\ref{rem:transl}, we will show $L\le_cL_0$.
Put $\pi(p,q)=(p\to q)\to q$.
Using the notation from \th\ref{def:fl}, \th\ref{rem:transl}, and from
the proof of \th\ref{thm:cpc-m}, we put $f(\alpha_n)=\beta_n$, where we define
inductively
$$
\beta_n:=(q\to q)\land
  \ET_{\vec\alpha_X,\vec\alpha_Y,\alpha_n\vdash_L\alpha_k}
   ((\vec\beta_Y\to\beta_k)\ot\vec\beta_X)
  \land\biggl(\Bigl(\ET_{\vec\alpha_Z\vdash_L\alpha_n}
         (\pi(p_n,q)\ot\vec\beta_Z)\Bigr)
     \to\pi(p_n,q)\biggr).$$
The first big conjunction in $\beta_n$ is taken over all $k<n$ and all
\emph{repetition-free} disjoint sequences $X$ and $Y$ consisting of
elements $i<n$ such that 
$\vec\alpha_X,\vec\alpha_Y,\alpha_n\vdash_L\alpha_k$, and similarly for the
second conjunction. (Here, if $X=\p{i_1,\dots,i_m}$, we define
$\vec\alpha_X$ to be the sequence
$\p{\alpha_{i_1},\dots,\alpha_{i_m}}$, and similarly for
$\vec\beta_X$.) If there are no $Z\sset n$ such that
$\vec\alpha_Z\vdash_L\alpha_n$, then the last conjunct of $\beta_n$ is
understood to be just $\pi(p_n,q)$.

Since $\beta_n(q/\bot)$ is classically equivalent to \eqref{eq:tr}, we
obtain immediately
$$f(\Gamma)\vdash_{L_0}f(\fii)\ \TO\ 
  f(\Gamma)\vdash_\CPC f(\fii)\ \TO\ \Gamma\vdash_L\fii.$$
In order to show
$$\Gamma\vdash_L\fii\ \TO\ f(\Gamma)\vdash_\FL{}f(\fii)
  \ \TO\ f(\Gamma)\vdash_{L_0}f(\fii),$$
it suffices to prove by
induction on $n$ that for every $k<n$ and every sequence $Z$ of
elements of $n$,
\begin{equation}\label{eq:6}
\vec\alpha_Z\vdash_L\alpha_k\ \TO\ {}\vdash_\FL{}\vec\beta_Z\to\beta_k
\end{equation}
(then $\vec\beta_Z\vdash_\FL{}\beta_k$ by modus ponens).
The statement is vacuously true for $n=0$. Assume that it holds for
$n$, we will prove it for $n+1$.
\begin{Cl}
Let $\p{L,\land,\lor,\cdot,\to,\ot,1}$ be a residuated lattice, and
$u\in L$.
\begin{enumerate}
\item\label{item:lu}
$L_u:=\{a\in L:au,ua\le a\}$ is closed under $\to,\ot,\land$ (as
well as $\cdot,\lor$, but we will not need this).
\item\label{item:qq}
If $u=q\to q$ for some $q\in L$, then $1\le u$, $u^2\le u$, and
$L_u$ contains $u$ as well as all elements of the form $\pi(a,q)$.
\item\label{item:bw}
$\FL{}$ proves $\beta_i\cdot\beta_j\to\beta_i$,
$\beta_j\cdot\beta_i\to\beta_i$, and
$\vec\beta_X\to(q\to q)$.
\end{enumerate}
\end{Cl}
\begin{Pf*}
\eqref{item:lu}:
Let $a,b\in L_u$.

We have $u(a\land b)\le ua\le a$ and $u(a\land b)\le ub\le b$, hence
$u(a\land b)\le a\land b$. The proof of $(a\land b)u\le a\land b$ is
symmetric.

Since $a(a\to b)u\le bu\le b$, we have $(a\to b)u\le a\to b$.
Similarly, $au(a\to b)\le a(a\to b)\le b$, hence $u(a\to b)\le a\to
b$.

The case of $a\ot b$ is symmetric.

\eqref{item:qq}: $1\le u$ is clear, and $u^2\le u$ (which implies
$u\in L_u$) is a special case of
\begin{equation}\label{eq:7}
(a\to b)(b\to c)\le a\to c.
\end{equation}
Put $p=(a\to q)\to q$. We have $p(q\to q)\le p$ from \eqref{eq:7}.
Also, $(a\to q)(q\to q)p\le (a\to q)p\le q$, hence
$(q\to q)p\le p$.

\eqref{item:bw}:
Consider a valuation $v$ in a residuated lattice $L$, and put
$u=v(q)\to v(q)$. Notice that $\beta_i$ is ultimately constructed from
formulas of the form $\pi(p_k,q)$ and $q\to q$ by means of
$\to,\ot,\land$, thus $v(\beta_i)\in L_u$ by \eqref{item:lu} and
\eqref{item:qq}. Clearly, $v(\beta_j)\le u$, hence $v(\beta_i)v(\beta_j)\le
v(\beta_i)$ and $v(\beta_j)v(\beta_i)\le v(\beta_i)$. Finally,
$v(\prod\vec\beta_{X^{-1}})\le u^{\lh X}\le u$ by \eqref{item:qq}.
\end{Pf*}
It follows from the claim that it is enough to prove \eqref{eq:6} for
repetition-free sequences $Z$ not containing $k$. The only interesting
cases are those involving $n$: i.e., $Z=X\cat n\cat Y$ or $k=n$.

Assume that $\vec\alpha_X,\alpha_n,\vec\alpha_Y\vdash_L\alpha_k$. Then the
definition of $\beta_n$ ensures
$\vdash_\FL{}\beta_n\to((\vec\beta_Y\to\beta_k)\ot\vec\beta_X)$, hence
$\vdash_\FL{}\beta_n\cdot\prod\vec\beta_{X^{-1}}\to(\vec\beta_Y\to\beta_k)$,
which in turn gives
$\vdash_\FL{}\prod\vec\beta_{X^{-1}}\to(\beta_n\to(\vec\beta_Y\to\beta_k))$
and $\vdash_\FL{}\vec\beta_X\to(\beta_n\to(\vec\beta_Y\to\beta_k))$.

Assume that $\vec\alpha_Z\vdash_L\alpha_n$. We have
$\vdash_\FL{}\vec\beta_Z\to(q\to q)$ by the Claim. Whenever
$\vec\alpha_X,\alpha_n,\vec\alpha_Y\vdash_L\alpha_k$, we have
$\vec\alpha_X,\vec\alpha_Z,\vec\alpha_Y\vdash_L\alpha_k$ by cut, hence
$\vdash_\FL{}\vec\beta_X\to(\vec\beta_Z\to(\vec\beta_Y\to\beta_k))$ by the
induction hypothesis. By a similar argument as above, this is equivalent to
$\vdash_\FL{}\vec\beta_Z\to((\vec\beta_Y\to\beta_k)\ot\vec\beta_X)$. Finally, that
$\vec\beta_Z$ implies the last conjunct of $\beta_n$ follows from
$\vdash_\FL{}\vec\beta_Z\to((\pi(p_n,q)\ot\vec\beta_Z)\to\pi(p_n,q))$.

\eqref{item:fle} follows immediately from \eqref{item:fl}, as
$(\fii\ot\psi)=(\psi\to\fii)$ in $\FL e$.

\eqref{item:bck}:
We define inductively
\begin{align*}
r_0&:=0,\\
r_{n+1}&:=1+n2^nr_n,\\
\ep_n&:=\biggl(\prod_{\ob\alpha_X,\alpha_n\vdash_L\alpha_k}(\ob\beta_X^{r_n}\to\beta_k)\biggr)
   \cdot\biggl(\Bigl(\prod_{\ob\alpha_Z\vdash_L\alpha_n}(\ob\beta_Z\to p_n)\Bigr)
               \to p_n\biggr),\\
\beta_n&:=(\ep_n\to q)\to q,
\end{align*}
where the products are taken over $X,Z\sset n$, $k<n$. (Unlike the
case of $\FL{}$, we can treat here
$X,Z$ as sets, because fusion is commutative.) 
We understand
$\ob\beta_X^{r_n}$ to be the multiset of formulas which contains $r_n$
copies of each formula $\beta_i$, $i\in X$ (again, the order does not
matter due to commutativity). That is, if $X=\{i_1,\dots,i_m\}$ (in an
arbitrary order), then $\ob\beta_X^{r_n}\to\beta_k$ stands for
$$\underbrace{\beta_{i_1}\to(\beta_{i_1}\cdots\to(\beta_{i_1}}_{r_n}\to
(\underbrace{\beta_{i_2}\cdots\to(\beta_{i_2}}_{r_n}\to
(\cdots\to
(\underbrace{\beta_{i_m}\cdots\to(\beta_{i_m}}_{r_n}\to\beta_k)))\cdots))).$$
Notice that 
fusion only appears in $\beta_n$ in premises of implications, hence
$\beta_n$ can be equivalently rewritten as a formula $f(\alpha_n)$
using only $\to$.

Since $\beta_n(q/\bot)$ is classically equivalent to \eqref{eq:tr}, we
have
$$f(\Gamma)\vdash_{L_0}f(\fii)\ \TO\ 
  f(\Gamma)\vdash_\CPC f(\fii)\ \TO\ \Gamma\vdash_L\fii.$$
By induction on $n$, we will show that
\begin{equation}\label{eq:8}
\ob\alpha_W\vdash_L\alpha_k\ \TO\ {}\vdash_\FL{ew}\ob\beta_W^{r_n}\to\beta_k
\end{equation}
holds for every $k<n$ and every $W\sset n$, which implies
$$\Gamma\vdash_L\fii\ \TO\ f(\Gamma)\vdash_\FL{ew}f(\fii)
  \ \TO\ f(\Gamma)\vdash_{L_0}f(\fii).$$
The statement is vacuously
true for $n=0$. Assume that it holds for $n$, we will prove it for
$n+1$. Since we have weakening, it suffices to consider the cases
$k=n$, $W\sset n$ and $k<n$, $W=X\cup\{n\}$, $X\sset n$.

Assume $\ob\alpha_X,\alpha_n\vdash_L\alpha_k$. Using the definition and
commutativity, we have
$\vdash_\FL{ew}\ob\beta_X^{r_n}\to(\ep_n\to\beta_k)$, which implies
$$\vdash_\FL{ew}\ob\beta_X^{r_n}\to((\beta_k\to q)\to(\ep_n\to q)).$$
The definition of $\beta_k$ gives $\vdash_\FL{ew}(\ep_k\to
q)\to(\beta_k\to q)$ using commutativity, hence
$$\vdash_\FL{ew}\ob\beta_X^{r_n}\to((\ep_k\to q)\to(\ep_n\to q)).$$
This implies $\vdash_\FL{ew}\ob\beta_X^{r_n}\to(((\ep_n\to q)\to
q)\to((\ep_k\to q)\to q))$, i.e.,
$\vdash_\FL{ew}\ob\beta_X^{r_n}\to(\beta_n\to\beta_k)$. We obtain
$$\vdash_\FL{ew}\ob\beta_X^{r_{n+1}}\to(\beta_n^{r_{n+1}}\to\beta_k)$$
by weakening, using $r_n\le r_{n+1}$.

Assume $\ob\alpha_W\vdash_L\alpha_n$. We have
$\vdash_\FL{ew}\ob\beta_W\to((\ob\beta_W\to p_n)\to p_n)$, hence
$$\vdash_\FL{ew}\ob\beta_W\to
  \biggl(\Bigl(\prod_{\ob\alpha_Z\vdash_L\alpha_n}(\ob\beta_Z\to p_n)\Bigr)\to p_n\biggr)$$
by weakening. Whenever $\ob\alpha_X,\alpha_n\vdash_L\alpha_k$, we have
$\ob\alpha_X,\ob\alpha_W\vdash_L\alpha_k$ by cut, hence
$$\vdash_\FL{ew}\ob\beta_W^{r_n}\to(\ob\beta_X^{r_n}\to\beta_k)$$
by the induction hypothesis and weakening. Since there are at most
$n2^n$ pairs $\p{X,k}$ such that $X\sset n$, $k<n$,
and $\ob\alpha_X,\alpha_n\vdash_L\alpha_k$, we have
$$\vdash_\FL{ew}\ob\beta_W^{n2^nr_n}\to
  \prod_{\ob\alpha_X,\alpha_n\vdash_L\alpha_k}(\ob\beta_X^{r_n}\to\beta_k).$$
Putting the pieces together, we have
$\vdash_\FL{ew}\ob\beta_W^{r_{n+1}}\to\ep_n$, hence
$$\vdash_\FL{ew}\ob\beta_W^{r_{n+1}}\to\beta_n.\qedhere$$
\end{Pf}
\begin{Rem}
Every consistent substitution-invariant extension of $\lgc{BCK}$ (in
the same language) is contained in $\CPC\res_\to$. This is no longer
true for $\FL e\res_{\to,\land}$, nevertheless one can modify the
proof above to show that \eqref{item:fl} and \eqref{item:fle} of
\th\ref{thm:univ} remain true when $\CPC$ is replaced with any
consistent substitution-invariant extension of
$\FL{}\res_{\to,\ot,\land}$ or $\FL e\res_{\to,\land}$, respectively.
We omit the details.
\end{Rem}

\section{Conclusion}\label{sec:conclusion}
Our results (\th\ref{thm:cpc-1,thm:univ,exm:kleene}) show that any countable
finitary deductive system can be conservatively translated into (among others):
\begin{itemize}
\item Classical propositional logic.
\item Intuitionistic, minimal, and intermediate logics.
\item Modal logics (classical or intuitionistic), including
variants such as temporal or epistemic logics.
\item Substructural logics, such as various extensions of $\FL{}$ or
linear logic.
\item Fuzzy and many-valued logics, such as $\lgc{MTL}$, $\lgc{BL}$
and their extensions (e.g., \L ukasiewicz logic).
\item Relevant logics, such as $\lgc R$.
\item Kleene's logic.
\item First-order (or higher-order) extensions of the above logics.
\item Implication fragments of many of the above logics.
\end{itemize}
This includes most of logical systems (fitting into the framework of
Tarski-style consequence relations) studied in the literature on
non-classical logic. We have also discovered some counterexamples,
namely $\CPC$ cannot be conservatively translated into its monotone or
affine fragments, or into the paraconsistent logic $\lgc{LP}$.

While there are still some loose ends left (most importantly, we were
unable to determine whether the logic $\lgc{BCI}=\FL e\res_\to$ is
universal, though it seems plausible), these results show that the
mere existence of a conservative translation of one logic into
another without further restrictions does not provide useful
information on the relationship of the two logics, and a more refined
criterion is needed to formalize the intuitive notion of translatability.

\subsection*{Acknowledgements}
I would like to thank Petr Cintula for various helpful suggestions on
a preliminary version of this paper, and the anonymous referees for
useful comments.

\bibliographystyle{mybib}
\bibliography{mybib}

\providecommand\gobble[1]{} {\catcode`\/=13
  \gdef/{\string/\futurelet\nexttoken\finishslash}
  \gdef\finishslash{\ifx\nexttoken/\else\penalty\relpenalty\fi}}
  \providecommand\url{\begingroup\catcode`\~=12 \catcode`\/=13 \finishurl}
  \def\finishurl#1{\texttt{#1}\endgroup}
  \providecommand\dotminus{\mathbin{\setbox0\hbox{$-$}\setbox2\hbox
  to\wd0{\hss$^{\mkern1mu\cdot}$\hss}\wd2=0pt\box2\box0}}
\providecommand{\bysame}{\leavevmode\hbox to5em{\hrulefill}\thinspace}
\providecommand\bibliographyhook{}
\begin{thebibliography}{10}
\bibliographyhook

\bibitem{cco}
Walter~A. Carnielli, Marcelo~E. Coniglio, and Itala M.~Loffredo D'Ottaviano,
  \emph{New dimensions on translations between logics}, Logica Universalis 3
  (2009), no.~1, pp.~1--18.

\bibitem{czel}
Janusz Czelakowski, \emph{Protoalgebraic logic}, Trends in Logic vol.~10,
  Kluwer, 2001.

\bibitem{dott-feit:mv}
Itala M.~Loffredo D'Ottaviano and H{\'e}rcules~A. Feitosa, \emph{Many-valued
  logics and translations}, Journal of Applied Non-classical Logics 9 (1999),
  no.~1, pp.~121--140.

\bibitem{dott-feit:para}
\bysame, \emph{Paraconsistent logics and translations}, Synthese 125 (2000),
  no.~1--2, pp.~77--95.

\bibitem{dott-feit:luk}
\bysame, \emph{Translating from {{\L}ukasiewicz}'s logics into classical logic:
  is it possible?}, Poznan Studies in the Philosophy of the Sciences and the
  Humanities 91 (2006), no.~1, pp.~157--168.

\bibitem{dott-feit:persp}
\bysame, \emph{Deductive systems and translations}, in: Perspectives on
  Universal Logic (J.-Y. B{\'e}ziau and A.~Costa-Leite, eds.), Polimetrica,
  2007, pp.~125--157.

\bibitem{feit-dott}
H{\'e}rcules~A. Feitosa and Itala M.~Loffredo D'Ottaviano, \emph{Conservative
  translations}, Annals of Pure and Applied Logic 108 (2001), pp.~205--227.

\bibitem{reslat}
Nikolaos Galatos, Peter Jipsen, Tomasz Kowalski, and Hiroakira Ono,
  \emph{Residuated lattices: An algebraic glimpse at substructural logics},
  Studies in Logic and the Foundations of Mathematics vol. 151, Elsevier,
  Amsterdam, 2007.

\bibitem{lau}
Dietlinde Lau, \emph{Function algebras on finite sets: A basic course on
  many-valued logic and clone theory}, Springer, New York, 2006.

\bibitem{mdt}
Till Mossakowski, R{\u a}zvan Diaconescu, and Andrzej Tarlecki, \emph{What is a
  logic translation?}, Logica Universalis 3 (2009), no.~1, pp.~95--124.

\bibitem{post}
Emil~L. Post, \emph{The two-valued iterative systems of mathematical logic},
  Annals of Mathematics Studies no.~5, Princeton University Press, Princeton,
  1941.

\bibitem{sh-sm}
D.~J. Shoesmith and Timothy~J. Smiley, \emph{Multiple-conclusion logic},
  Cambridge University Press, 1978.

\bibitem{sos}
Jairo~J. da~Silva, Itala M.~Loffredo D'Ottaviano, and Ant{\^o}nio~M. Sette,
  \emph{Translations between logics}, in: Models, Algebras, and Proofs
  (X.~Caicedo and C.~Montenegro, eds.), Lecture Notes in Pure and Applied
  Mathematics vol. 203, Marcel Dekker, New York, 1999, pp.~435--448.

\end{thebibliography}
\end{document}